# q-Abel polynomials

*Johann Cigler*

Fakultät für Mathematik
Universität Wien
A-1090 Wien, Nordbergstraße 15

Johann.Cigler@univie.ac.at

**Abstract**

This note gives a simple approach to $q-$ analogues of some results associated with Abel polynomials.

## 0. Introduction

In this note I want to give a simple approach to $q-$ analogues of the following well-known results about Abel polynomials:
Let $a_n(x,a) = x(x-na)^{n-1}$ be the Abel polynomials.
N.H. Abel [1] has found the beautiful formula

$$(x+y)^n = \sum_{k=0}^{n} \binom{n}{k} a_k(x,a)(y+ka)^{n-k}. \tag{0.1}$$

We want to state it in a slightly more general form by changing $x \to x-b$, $y \to y+b$ and defining

$$a_n(x,a,b) = a_n(x-b,a) = (x-b)(x-b-na)^{n-1}. \tag{0.2}$$

Then we get

$$(x+y)^n = \sum_{k=0}^{n} \binom{n}{k} a_k(x,a,b)(y+ka+b)^{n-k}. \tag{0.3}$$

If we denote by $\partial$ the differentiation operator we see that

$$a_n(x,a,b) = (1+a\partial)(x-b-na)^n. \tag{0.4}$$

By Taylor's theorem which may be stated as $p(y+x) = e^{x\partial_y} p(y)$ for polynomials $p(y)$ formula (0.3) is equivalent with

$$e^{x\partial_y} p(y) = \sum_k \frac{a_k(x,a,b)}{k!} e^{(b+ka)\partial_y} \left(\partial_y\right)^k p(y) \tag{0.5}$$

for all polynomials $p(y)$. This gives the operator identity

$$e^{x\partial_y} = \sum_k \frac{a_k(x,a,b)}{k!} e^{(b+ka)\partial_y} \left(\partial_y\right)^k \tag{0.6}$$

which by the isomorphism $\partial \leftrightarrow z$ is equivalent with the identity of formal power series

$$e^{xz} = \sum_k \frac{a_k(x,a,b)}{k!} z^k e^{(b+ka)z}. \tag{0.7}$$

The linear operator

$$Q = \partial e^{a\partial}, \tag{0.8}$$

called Abel operator, satisfies

$$Q a_n(x,a,b) = n a_{n-1}(x,a,b). \tag{0.9}$$

Since $a_n(0,a) = [n=0]$ we get

$$LQ^k a_n(x,a) = L\partial^k e^{ka\partial} a_n(x,a) = n![k=n], \tag{0.10}$$

where $L$ denotes the linear functional on the polynomials $p(x)$ defined by $Lp = p(0)$ and $[k=n]$ is Knuth's symbol defined by $[k=n]=1$ if $k=n$ is true and $[k=n]=0$ if $k \neq n$.

This implies that the coefficients of the expansion

$$f(z) = \sum_k \frac{c_k}{k!} z^k e^{kaz} \tag{0.11}$$

are given by the Lagrange formula

$$c_n = L\partial^{n-1} e^{-nax} f'(x). \tag{0.12}$$

In the same way we get that the coefficients of the expansion

$$\frac{f(z)}{1+az} = \sum_k \frac{c_k}{k!} z^k e^{kaz} \tag{0.13}$$

are given by the formula of Lagrange-Bürmann

$$c_n = L\partial^n e^{-nax} f(x). \tag{0.14}$$

A special case is

$$\frac{e^{xz}}{1-az} = \sum_{k \geq 0} \frac{(x+ak)^k}{k!} \left(ze^{-az}\right)^k. \tag{0.15}$$

This can also be written in the form

$$\frac{1}{1-az} = \sum_{k \geq 0} \frac{(x+ak)^k}{k!} z^k e^{-(ak+x)z}. \tag{0.16}$$

Expanding $e^{-(ak+x)z}$ into a power series and comparing coefficients of $z^n$ shows that this is equivalent with

$$n! a^n = \sum_{k=0}^n (-1)^{n-k} \binom{n}{k} (x+ak)^n. \tag{0.17}$$

Historical notes about some of these formulas can be found in [10].

The first $q$-analogues of some of these results have been given by F. H. Jackson ([5]) in 1910. In [3] I have given another proof depending on Rota's Finite Operator Calculus ([8]). The case $b \neq 0$ has first been considered by J. Hofbauer in [4]. It also appeared in papers by W.P. Johnson ([6]), B. Bhatnagar and St.C. Milne ([2]), C. Krattenthaler and M. Schlosser ([7]) and M. Schlosser ([9]). In the following I want to give a self-contained exposition and some simplifications of these results. I want to thank Michael Schlosser for some useful comments.

## 1. q-Abel polynomials

Let $D$ or $D_x$ be the $q$-differentiation operator, defined by $Df(x) = \dfrac{f(x) - f(qx)}{(1-q)x}$. Instead of $Df(x)$ we also write $f'(x)$.

In order to simplify notation we set $(y \dotplus x)^n := \prod_{j=0}^{n-1}(y + q^j x)$ and $(y \dotminus x)^n := \prod_{j=0}^{n-1}(y - q^j x)$.

The other notations from $q$-calculus are the usual ones. The $q$-binomial coefficients are denoted by $\begin{bmatrix} n \\ k \end{bmatrix} = \dfrac{[n]!}{[k]![n-k]!}$ for $0 \leq k \leq n$. Here $[n] = \dfrac{1-q^n}{1-q}$ and $[n]! = \prod_{j=1}^{n}[j]$.

We need the analogues of the exponential series $e(z) = \sum_{k \geq 0} \dfrac{z^k}{[k]!}$ and $E(z) = \sum_{k \geq 0} q^{\binom{k}{2}} \dfrac{z^k}{[k]!}$, which are related by $e(z)E(-z) = 1$ and the well-known facts that $\sum_{k \geq 0} \dfrac{(x \dotminus y)^k}{[k]!} z^k = \dfrac{e(xz)}{e(yz)}$ and $E(aD)y^n = (y+a)\cdots(y+q^{n-1}a)$.

The following $q$-analogue of (0.3) holds:

$$\prod_{j=0}^{n-1}(y + q^j x) = \sum_{k=0}^{n} \begin{bmatrix} n \\ k \end{bmatrix} A_k(x,a,b) \prod_{j=0}^{n-k-1}(y + q^j[k]a + q^{k+j}b). \tag{1.1}$$

Here

$$A_n(x,a,b) = (x-b)\prod_{j=1}^{n-1}\left(q^j x - a[n] - q^n b\right) \tag{1.2}$$

is a $q$-analogue of the general Abel polynomial $a_n(x,a,b)$.

We will write (1.1) in the form

$$(y \dotplus x)^n = \sum_{k=0}^{n} \begin{bmatrix} n \\ k \end{bmatrix} A_k(x,a,b)(y \dotplus ([k]a + q^k b))^{n-k} \tag{1.3}$$

or

$$E(xD_y)y^n = \sum_k \dfrac{A_k(x,a,b)}{[k]!} E\left(\left([k]a + q^k b\right)D_y\right) D_y^k y^n. \tag{1.4}$$

By the isomorphism $D \leftrightarrow z$ it is also equivalent with the identity of formal power series

$$E(xz) = \sum_k \frac{A_k(x,a,b)}{[k]!} E(([k]a + q^k b)z)z^k. \tag{1.5}$$

For $b = 0$ this $q$ – Abel theorem has been found by F. H. Jackson ([5]). The general case has been considered by J. Hofbauer ([4]), W.P. Johnson ([6]) and is also contained in more general results by B. Bhatnagar and St.C. Milne ([2]) and by C. Krattenthaler and M. Schlosser ([7]) and M. Schlosser ([9]).

In order to prove (1.2) we write (1.3) with unknown coefficients $A_k(x,a,b)$ and try to determine their values.
We consider this formula for $n \to n-1$ and multiply with $c[n]$ for some constant $c$.
We thus get

$$c[n](y \dotplus x)^{n-1} = c \sum_{k=0}^{n} \begin{bmatrix} n-1 \\ k \end{bmatrix} [n] A_k(x,a,b)(y \dotplus ([k]a + q^k b))^{n-k-1}$$

$$= \sum_{k=0}^{n} \begin{bmatrix} n-1 \\ n-1-k \end{bmatrix} \frac{[n]}{[n-k]} A_k(x,a,b)(y \dotplus ([k]a + q^k b))^{n-k} \frac{c[n-k]}{y + q^{n-k-1}([k]a + q^k b)}.$$

Comparing with (1.3) we see that the first $n-1$ terms coincide if $\frac{c[n-k]}{y + q^{n-k-1}([k]a + q^k b)} = 1$.

This gives $y = c[n-k] - q^{n-1-k} a[k] - q^{n-1}b = \frac{c(1-q^{n-k}) - q^{n-1-k}(1-q^k)a}{1-q} - q^{n-1}b$.

This is independent on $k$ if $q^{n-k}c + q^{n-1-k}a = 0$, i.e. $c = -\frac{a}{q}$.

We then get

$$y = -\frac{a}{q} \frac{(1-q^{n-k}) + q^{n-k}(1-q^k)}{1-q} - q^{n-1}b = -\frac{a}{q}[n] - q^{n-1}b.$$

Therefore

$$A_n(x,a,b) = q^{n-1}(x-b) \prod_{j=0}^{n-2} \left( q^j x - \frac{a}{q}[n] - q^{n-1}b \right) = (x-b) \prod_{j=1}^{n-1} \left( q^j x - a[n] - q^n b \right).$$

For $b = 0$ we get

$$A_n(x,a,0) = x \prod_{j=1}^{n-1} \left( q^j x - [n]a \right) \tag{1.6}$$

for $n > 0$ and $A_0(x,a,q) = 1$.

This implies Jackson's identity.

By letting $a \to a + (1-q)b$ we get

$$G_n(x,a,b) = A_n(x, a+(1-q)b, b) = (x-b)\prod_{j=1}^{n-1}(q^j x - [n]a - b). \tag{1.7}$$

This is another form of the general $q$ – Abel polynomials. Whereas $A_n(x,a,b)$ and $G_n(x,a,b)$ are equivalent some formulas become simpler by using $G_n(x,a,b)$.

The polynomials $G_n(x,a,b)$ have first been considered by J. Hofbauer ([4]). They are also the special case $h = 0$ of the polynomials $a_n(x;b,h,w,q)$, which have been studied by W.P. Johnson in [6].

The corresponding identity is

$$\prod_{j=0}^{n-1}(y+q^j x) = \sum_{k=0}^{n}\begin{bmatrix}n\\k\end{bmatrix}G_k(x,a,b)\prod_{j=0}^{n-k-1}(y+q^j[k]a+q^j b). \tag{1.8}$$

This is equivalent with identity (8.4) in [7]. There it is stated in the form

$$(c;q)_n = \sum_{k=0}^{n}\begin{bmatrix}n\\k\end{bmatrix}(-1)^k q^{\binom{k}{2}}c^k \frac{1-(a+b)}{1-(q^{-k}a+b)}\left(q^{-k}a+b;q\right)_k \left(c(q^k a + b);q\right)_{n-k},$$

where as usual $(x;q)_n = (1-x)(1-qx)\cdots(1-q^{n-1}x)$. To obtain (1.8) make the substitutions
$$a \to \frac{b(1-q)-a}{(1-q)x}, b \to \frac{a}{(1-q)x}, c \to -\frac{x}{y}.$$

The corresponding formula
$$E(xz) = \sum_k \frac{G_k(x,a,b)}{[k]!} E(([k]a+b)z)z^k \tag{1.9}$$

has been obtained in [4] and is also equivalent with formula [7], (7.4).

## 2. Abel expansions

If we apply the operator $D$ to (1.9) we get

$$zE(qxz) = DE(xz) = \sum_{n\geq 1}\frac{DG_n(x,a,b)}{[n]!}z^n E((b+[n]a)z)$$

or

$$E(qxz) = \sum_{n\geq 0}\frac{DG_{n+1}(x,a,b)}{[n+1]!}z^n E((b+[n+1]a)z)).$$

On the other hand (1.9) also gives

$$E(qxz) = \sum_{n\geq 0} \frac{G_n(qx, qa, b+a)}{[n]!} z^n E((b+[n+1]a)z)).$$

Comparing these two formulas we get

$$DG_{n+1}(x,a,b) = [n+1]G_n(qx, qa, b+a).$$

By induction this implies

$$D^k G_n(x,a,b) = q^{\binom{k}{2}} \frac{[n]!}{[n-k]!} G_{n-k}(q^k x, q^k a, b+[k]a). \qquad (2.1)$$

This means that

$$D^n G_n(x,a,b) = q^{\binom{n}{2}} [n]!$$

and for $k < n$

$$D^k G_n(x,a,b) = q^{\binom{k}{2}} \frac{[n]!}{[n-k]!} (q^k x - [k]a - b) \prod_{j=k+1}^{n-1} (q^j x - [n]a - b).$$

An important consequence is

$$G_n^{(k)}\left(\frac{b+[k]a}{q^k}, a, b\right) = q^{\binom{k}{2}} [k]! [k=n]. \qquad (2.2)$$

Thus each polynomial $f(x)$ has the following Abel expansion

$$f(x) = \sum_k \frac{f^{(k)}(q^{-k}(b+[k]a))}{[k]!} q^{-\binom{k}{2}} G_k(x,a,b). \qquad (2.3)$$

Choosing $f(x) = \prod_{j=0}^{n-1}(y + q^j x)$ we get again (1.8).

By choosing $f(x) = G_n(x, -a, -y-b)$ we get

$$G_n(x,-a,-y-b) = \sum_{k=0}^n \begin{bmatrix} n \\ k \end{bmatrix} G_k(x,-a,-b) y \prod_{j=1}^{n-k-1} \left(y + (1-q^j)b + ([n]-q^j[k])a\right). \qquad (2.4)$$

This is the special case $h=0$ of Theorem 4 by W.P. Johnson ([6]).

The expansion (2.3) also holds for formal power series. If we choose $f(x) = E(xz)$ we get (1.9).

From (1.9) we get e.g.

$$z^n = \sum_{k\geq 0} \frac{(-1)^k}{[k]!} ([n]a+b)([n+k]a+b)^{k-1} z^{n+k} E(([n+k]a+b)z)$$

by applying $q$-differentiation $k$ times and setting $x=0$.

## 3. A q - analogue of the Abel operator

Let
$$w_n(x,a,b) = \prod_{j=0}^{n-1}\left(q^j x - [n]a - b\right). \tag{3.1}$$

This can be written in the form

$$w_n(x,a,b) = \sum_{k=0}^{n}(-1)^k \begin{bmatrix} n \\ k \end{bmatrix} q^{\binom{n-k}{2}} ([n]a+b)^k x^{n-k} = q^{\binom{n}{2}} \sum_{k=0}^{n}(-1)^k \begin{bmatrix} n \\ k \end{bmatrix} q^{\binom{k}{2}} \left(\frac{[n]a+b}{q^{n-1}}\right)^k x^{n-k}$$

$$= q^{\binom{n}{2}} E\left(-\frac{[n]a+b}{q^{n-1}} D\right) x^n.$$

In analogy to (0.4) we get

$$G_n(x,a,b) = (1+aD)\prod_{j=0}^{n-1}\left(q^j x - [n]a - b\right) = (1+aD)w_n(x,a,b). \tag{3.2}$$

Therefore we have

$$G_n(x,a,b) = (1+aD)w_n(x,a,b) = (1+aD)q^{\binom{n}{2}} E\left(-\frac{[n]a+b}{q^{n-1}} D\right) x^n$$

$$= q^{\binom{n}{2}} E\left(-\frac{[n]a+b}{q^{n-1}} D\right)\left(x^n + [n]ax^{n-1}\right).$$

If we write $S_n(x,a) = q^{-\binom{n}{2}} e\left(\frac{[n]a+b}{q^{n-1}} D\right) G_n(x,a,b) = x^n + [n]ax^{n-1}$, then

$DS_n(x,a) = [n]S_{n-1}(x,a)$.

Let now $Q_n$ be the operator

$$Q_n = \frac{D}{q^{n-1}} \frac{e\left(\frac{[n]a+b}{q^{n-1}} D\right)}{e\left(\frac{[n-1]a+b}{q^{n-2}} D\right)}. \tag{3.3}$$

Then we get as a $q$ − analogue of (0.9)

$$Q_n G_n(x,a,b) = [n]G_{n-1}(x,a,b). \tag{3.4}$$

Therefore the operators $Q_n$ can be interpreted as $q$ − analogues of the Abel operator $Q$. Unfortunately they are depending on $n$.

(3.4) is a consequence of

$$\frac{D}{q^{n-1}} \frac{e\left(\frac{[n]a+b}{q^{n-1}}D\right)}{e\left(\frac{[n-1]a+b}{q^{n-2}}D\right)} G_n(x,a,b) = \frac{1}{q^{n-1}} E\left(-\frac{[n-1]a+b}{q^{n-2}}D\right) q^{\binom{n}{2}} DS_n(x,a,b)$$

$$= q^{\binom{n-1}{2}} E\left(-\frac{[n-1]a+b}{q^{n-2}}D\right) [n]S_{n-1}(x,a,b) = [n]G_{n-1}(x,a,b).$$

The operator $Q_n$ can also be written as

$$Q_n = \sum_{k \geq 0} \frac{\prod_{j=0}^{k-1}\left(([j+1]-q^n[j])a + (1-q^{j+1})b\right)}{[k]!} \left(\frac{D}{q^{n-1}}\right)^{k+1}. \qquad (3.5)$$

For

$$[j+1] - q^n[j] = \frac{1-q^{j+1} - q^n(1-q^j)}{1-q} = \frac{(1-q^n) - q^{j+1}(1-q^{n-1})}{1-q} = [n] - q^{j+1}[n-1]$$

and therefore

$$([j+1]-q^n[j])a + (1-q^{j+1})b = ([n]-q^{j+1}[n-1])a + (1-q^{j+1})b = ([n]a+b) - q^{j+1}([n-1]a+b).$$

This implies

$$\prod_{j=0}^{k-1}\left(([j+1]-q^n[j])a + (1-q^{j+1})b\right) = \left(([n]a+b) \div q([n-1]a+b)\right)^k.$$

## 4. A q - Lagrange formula

For a formal power series $f(z)$ we want to find the coefficients in the expansion

$$f(z) = \sum_n \frac{c_n}{[n]!} z^n E([n]az). \qquad (4.1)$$

Consider first the expansion of $f(z) = e(xz)$. Let

$$e(xz) = \sum_k \frac{B_k(x,a)}{[k]!} z^k E([k]az). \qquad (4.2)$$

We know from (1.9) that $E(xz) = \sum_k \frac{A_k(x,a,0)}{[k]!} E([k]az) z^k$.

If we let $V$ be the linear operator defined by $V q^{\binom{n}{2}} x^n = x^n$, then it is clear that

$$B_n(x,a) = V A_n(x,a,0) = V \sum_{k=0}^{n-1} \begin{bmatrix} n-1 \\ k \end{bmatrix} (-1)^k ([n]a)^k q^{\binom{n-k}{2}} x^{n-k} = \sum_{k=0}^{n-1} \begin{bmatrix} n-1 \\ k \end{bmatrix} (-1)^k ([n]a)^k x^{n-k}.$$

Therefore we get

$$B_n(x,a) = x \sum_{j=0}^{n-1} (-1)^j \begin{bmatrix} n-1 \\ j \end{bmatrix} a^j x^{n-1-j} [n]^j = x e(-[n]aD) x^{n-1}. \tag{4.3}$$

Here we have a direct analogue of the formula $a_n(x,a) = x e^{-na\partial} x^{n-1}$.

A possible disadvantage is that there is no simple factorization.

The property we are interested in is

$$L E([k]aD) D^k B_n(x,a) = [n]! [k = n], \tag{4.4}$$

which generalizes (0.10).

To prove it observe that $D^k e(xz) = z^k e(xz)$. Therefore

$$E([k]aD) D^k e(xz) = \sum_k \frac{E([k]aD) D^k B_k(x,a)}{[k]!} z^k E([k]az)$$

and for $x = 0$

$$E([k]az) z^k = \sum_k \frac{L E([k]aD) D^k B_k(x,a)}{[k]!} z^k E([k]az).$$

Comparing coefficients we get (4.4).

This implies the following

$q$ – **Lagrange formula** ([3]):

*The coefficients $c_n$ in the expansion*

$$f(x) = \sum_n \frac{c_n}{[n]!} x^n E([n]ax) \tag{4.5}$$

*are given by*

$$c_n = L f'(D) e(-[n]aD) x^{n-1} = L D^{n-1} e(-[n]ax) f'(x). \tag{4.6}$$

For by (4.4) we have $c_n = L f(D) B_n(x,a) = L f(D) x e(-[n]aD) x^{n-1} = L D^{n-1} e(-[n]aD) f'(x)$.

The last equation follows from $L(D^k x^n) = [n]![k=n] = L(D^n x^k)$ and the $q$-Pincherle derivative

$$f(D)\underline{x} - \underline{x}f(qD) = f'(D). \tag{4.7}$$

Here $\underline{x}$ denotes the operator multiplication by $x$.

This well-known fact follows from $(D\underline{x} - q\underline{x}D)x^n = Dx^{n+1} - qxDx^n = ([n+1]-q[n])x^n = x^n$, which implies $D\underline{x} - \underline{x}qD = 1$ by induction $D^n\underline{x} - \underline{x}q^nD^n = nD^{n-1}$.

*More generally let $B_n(x,a,b) = VA_n(x,a,b)$. Then the coefficients of the expansion*

$$f(z) = \sum_n \frac{c_n}{[n]!} z^n E\left(([n]a + q^n b)z\right) \tag{4.8}$$

*are given by*

$$c_n = Lf(D)B_n(x,a,b) = Lf(D)xe\left(-([n]a+q^n b)D\right)x^{n-1} - q^{n-1}bLf(D)e\left(-\left(\frac{[n]a+q^n b}{q}\right)D\right)x^{n-1}$$

$$= LD^{n-1}e\left(-([n]a+q^n b)x\right)f'(x) - q^{n-1}bLD^{n-1}e\left(-\left(\frac{[n]a+q^n b}{q}\right)x\right)f(x).$$

For the proof observe that

$$LE\left((q^k b + [k]a)D\right)D^k B_n(x,a,b) = [n]![k=n] \tag{4.9}$$

by the same argument as above and that

$$B_n(x,a,b) = VA_n(x,a,b) = V(x-b)\sum_{k=0}^{n-1}\begin{bmatrix}n-1\\k\end{bmatrix}(-1)^k ([n]a+q^n b)^k q^{\binom{n-k}{2}} x^{n-1-k}$$

$$= V\sum_{k=0}^{n-1}\begin{bmatrix}n-1\\k\end{bmatrix}(-1)^k ([n]a+q^n b)^k q^{\binom{n-k}{2}} x^{n-k} - bV\sum_{k=0}^{n-1}\begin{bmatrix}n-1\\k\end{bmatrix}(-1)^k ([n]a+q^n b)^k q^{\binom{n-1-k}{2}} q^{n-1-k} x^{n-1-k}$$

$$= \sum_{k=0}^{n-1}\begin{bmatrix}n-1\\k\end{bmatrix}(-1)^k ([n]a+q^n b)^k x^{n-k} - q^{n-1}b\sum_{k=0}^{n-1}\begin{bmatrix}n-1\\k\end{bmatrix}(-1)^k \left(\frac{[n]a+q^n b}{q}\right)^k x^{n-1-k}$$

$$= xe\left(-([n]a+q^n b)D\right)x^{n-1} - q^{n-1}be\left(-\left(\frac{[n]a+q^n b}{q}\right)D\right)x^{n-1}.$$

$B_n(x,a,b)$ can also be written in the form

$$B_n(x,a,b) = \sum_{k=0}^n (-1)^k x^{n-k}\begin{bmatrix}n\\k\end{bmatrix}(q^n b + [n]a)^{k-1}(q^{n-k}b + [n-k]a).$$

In order to obtain an analogue of the **Lagrange-Bürmann formula** we note that

$$\left(1+\frac{a}{q}D\right)e\left(-\frac{q^n b+[n]a}{q}D\right)x^n = xe\left(-\left([n]a+q^n b\right)D\right)x^{n-1} - q^{n-1}be\left(-\left(\frac{[n]a+q^n b}{q}\right)D\right)x^{n-1} \quad (4.10)$$
$$= B_n(x,a,b).$$

This follows from the $q$ – Pincherle derivative

$$\left(e\left(-\frac{q^n b+[n]a}{q}D\right)\underline{x} - \underline{x}e(-(q^n b+[n]a)D)\right)x^{n-1} = -\frac{q^n b+[n]a}{q}e\left(-\frac{q^n b+[n]a}{q}D\right)x^{n-1}.$$

Thus we get a

**q - Lagrange - Bürmann type formula:**

*The coefficients of*

$$\frac{f(z)}{1+\frac{az}{q}} = \sum_k \frac{c_k}{[k]!}z^k E\left(\left(q^k b+[k]a\right)z\right) \quad (4.11)$$

*are given by*

$$c_n = LD^n e\left(-\frac{q^n b+[n]a}{q}x\right)f(x). \quad (4.12)$$

This is an immediate consequence of

$$c_n = Lf(D)\frac{1}{1+\frac{aD}{q}}B_n(x,a,b) = Lf(D)e\left(-\frac{q^n b+[n]a}{q}D\right)x^n = LD^n e\left(-\frac{q^n b+[n]a}{q}x\right)f(x).$$

If we choose $f(z) = E(-yz)$ we get

$$c_n = LD^n e\left(-\frac{q^n b+[n]a}{q}x\right)E(-xy) = LD^n \frac{e\left(-\frac{q^n b+[n]a}{q}x\right)}{e(xy)} = \left(-\frac{q^n b+[n]a}{q} \dot{-} y\right)^n$$

This is equivalent with

$$\frac{E(xz)}{1-az} = \sum_k \frac{(q^k b+[k]a \dot{+} x)^k}{[k]!} z^k E(-q(q^k b+[k]a)z). \quad (4.13)$$

This $q$ – analogue of (0.15) has been found by C. Krattenthaler and M. Schlosser ( cf. [7] and [9], (5.4)) in another context.

## 5. Other methods of proof

We give now another proof of formula (4.13) by using a $q$-analogue of the difference operator:

Let $U$ be the linear operator on the polynomials in $q^n$ defined by

$$Uq^{in} = q^{i(n-1)} \tag{5.1}$$

for all $i \in \mathbb{N}$. As a special case we get $U[n]^m = [n-1]^m$.

Define now

$$\Delta^k = (1-qU)\cdots(1-q^kU). \tag{5.2}$$

Then it is clear that

$$\Delta^k q^{in} = 0 \tag{5.3}$$

for $k \geq i > 0$ and

$$\Delta^k 1 = (1-q)^k [k]!. \tag{5.4}$$

Furthermore

$$\Delta^k [n]^m = [k]!(1-q)^{k-m} \tag{5.5}$$

for $k \geq m$ and

$$\Delta^k q^{in} [n]^m = 0 \tag{5.6}$$

for $k \geq m+i$ if $i > 0$.

It suffices to show (5.5), which follows from

$$\Delta^k [n]^m = \frac{1}{(1-q)^m} \Delta^k (1-q^n)^m = \frac{1}{(1-q)^m} \sum (-1)^j \binom{m}{j} \Delta^k q^{nj} = \frac{\Delta^k 1}{(1-q)^m} = \frac{(1-q)^k [k]!}{(1-q)^m}.$$

This again implies

$$\Delta^n \left( q^{nj}(q^n x + [n]a)^{n-j} \right) = 0 \tag{5.7}$$

if $j > 0$ and

$$\Delta^n \left( (q^n x + [n]a)^n \right) = [n]!a^n. \tag{5.8}$$

We conclude that

$$\sum_{k=0}^{n} \begin{bmatrix} n \\ k \end{bmatrix} (-1)^k \left( q^{n-k}b + [n-k]a \dotplus x \right)^{n-k} \left( x \dotplus q(q^{n-k}b + [n-k]a) \right)^k = [n]!a^n. \tag{5.9}$$

For

$$\sum_{k=0}^{n}\begin{bmatrix}n\\k\end{bmatrix}(-1)^k\left(q^{n-k}b+a[n-k]\dotplus x\right)^{n-k}\left(x\dotplus q(q^{n-k}b+a[n-k])\right)^k$$

$$=\sum_{k=0}^{n}\begin{bmatrix}n\\k\end{bmatrix}(-1)^k q^{\binom{k+1}{2}}\left(q^{n-k}b+a[n-k]\dotplus q^{-k}x\right)^n$$

$$=\sum_{k=0}^{n}\begin{bmatrix}n\\k\end{bmatrix}(-1)^k q^{\binom{k+1}{2}}\sum_{j=0}^{n}q^{-kj}x^j q^{\binom{j}{2}}\left(q^{n-k}b+a[n-k]\right)^{n-j}$$

$$=\sum_{j=0}^{n}q^{-nj}x^j q^{\binom{j}{2}}\Delta^n\left(q^{nj}\left(q^n b+a[n]\right)^{n-j}\right)=\Delta^n\left(\left(q^n b+a[n]\right)^n\right)=[n]!a^n.$$

By comparing coefficients we see that (5.9) is equivalent with

$$\frac{1}{1-az}=\sum_k\frac{(q^k b+[k]a\dotplus x)^k}{[k]!}z^k\sum_j\frac{\left(x\dotplus q(a[k]+q^k b)\right)^j}{[j]!}(-z)^j$$

$$=\sum_k\frac{(q^k b+[k]a\dotplus x)^k}{[k]!}z^k\frac{e(-xz)}{e(q(q^k b+[k]a)z)}.$$

This is equivalent with (4.13).

By applying the $q$-differentiation operator $k$ times and then setting $x=0$ we get from (4.13)

$$\frac{z^n}{1-az}=\sum_{k\geq 0}\frac{1}{[k]!}\left([n+k]a+q^{n+k}b\right)^k z^{n+k}E\left(-q\left([n+k]a+q^{n+k}b\right)z\right). \quad (5.10)$$

By considering the isomorphism $z\to D_y$ and applying it to $y^n$ (4.13) gives

$$(y\dotplus x)^n=\sum_{k=0}^{n}\begin{bmatrix}n\\k\end{bmatrix}\left(q^k b+[k]a\dotplus x\right)^k v(n,k,a,b,y), \quad (5.11)$$

with $v(n,k,a,b,y)=\left(y\dotminus q(q^k b+[k]a)\right)^{n-k-1}(y-q^n b-[n]a)$ for $k<n$ and $v(n,n,a,b,y)=1$.

If we substitute $a\to a+(1-q)b$ in (4.13) we get

$$\frac{E(xz)}{1-\left(a+(1-q)b\right)z}=\sum_k\frac{(b+[k]a\dotplus x)^k}{[k]!}z^k E(-q(b+[k]a)z). \quad (5.12)$$